\documentclass[12pt]{amsart}
\newtheorem{theorem}{Theorem}
\newtheorem{proposition}{Proposition}
\newtheorem{lemma}{Lemma}

\newcommand{\sm}{\setminus}

\newcommand{\F}{{\mathcal F}}

\newcommand{\C}{{\bf C}}
\newcommand{\R}{{\bf R}}
\newcommand{\Z}{{\bf Z}}

\newcommand{\CP}{{\bf CP}}
\newcommand{\HH}{{\mathcal H}}
\newcommand{\E}{{\mathcal E}}
\newcommand{\thm}{\begin{theorem}}
\newcommand{\etheorem}{\end{theorem}}
\newcommand{\prop}{\begin{proposition}}
\newcommand{\eprop}{\end{proposition}}

\author{V.~A.~Vassiliev}

\email{vva@mi.ras.ru}

\thanks{supported in part by INTAS
(project 96-0713) and RFBR (project 98-01-00555)}

\title{Spaces of Hermitian operators with simple spectra
and their finite-order cohomology}

\date{Revised version published in 2003}

\begin{document}

\begin{abstract}
V.~I.~Arnold \cite{hall} studied the topology of spaces of Hermitian operators
in $\C^n$ with non-simple spectra in a relation with the theory of adiabatic
connections and the quantum Hall effect. The natural filtration of these spaces
by the sets of operators with fixed numbers of eigenvalues defines the spectral
sequence, providing interesting combinatorial and homological information on
this stratification.

We construct a different spectral sequence, also counting the homology groups
of these spaces and based on the universal techniques of {\em topological order
complexes} and {\em conical resolutions} of algebraic varieties, generalizing
the combinatorial inclusion-exclusion formula and similar to the construction
of finite-order knot invariants.

This spectral sequence degenerates at the term $E_1$, is (conjecturally)
multiplicative, and as $n$ grows then it converges to a stable spectral
sequence counting the cohomology of the space of infinite Hermitian operators
without multiple eigenvalues, all whose terms $E^{p,q}_r$ are finitely
generated. It allows us to define the finite-order cohomology classes of this
space, and to apply the well-known facts and methods of the topological theory
of flag manifolds to the problems of geometrical combinatorics, especially
concerning the continuous partially ordered sets of subspaces and flags.
\end{abstract}

\maketitle

\section{Introduction}

Denote by $\HH(n)$ the space of Hermitian operators in $\C^n$; this is a
$n^2$-dimensional real vector space. The {\em discriminant variety} $\Sigma
\equiv \Sigma(n) \subset \HH(n)$ consists of operators with at least one
eigenvalue of multiplicity $\ge 2$; this is a subvariety of codimension 3. The
Alexander duality
\begin{equation}
\label{alex} \bar H_i (\Sigma(n)) \simeq \tilde H^{n^2-i-1} (\HH(n) \sm
\Sigma(n))
\end{equation}
relates its Borel--Moore homology groups (i.e. homology groups of locally
finite chains) with the standard cohomology groups of the complementary space
of matrices with simple spectra.

V.~I.~Arnold \cite{hall} studied these homology groups, considering the natural
filtration of the space $\Sigma(n)$ in accordance with the numbers of different
eigenvalues of operators. Although the answer is known (there is a group
isomorphism
\begin{equation}
\label{split} H^*(\HH(n)\sm \Sigma(n)) \cong H^*(\CP^{n-1} \times \CP^{n-2}
\times \cdots \times \CP^1)),
\end{equation}
and the ring structure also is easy\footnote{The reduced cohomology ring is
isomorphic to $\Z[c_1, \dots, c_n]/\mbox{Sym}$, where all $c_i$ are Chern
classes of bundles of proper subspaces, and $\mbox{Sym}$ is the ideal generated
by symmetric polynomials, see \cite{borel}}, this study gives interesting
information on the topology of the discriminant set. In particular, the
spectral sequence, defined by this filtration, degenerates at the second term
$E^2$ (see \cite{ShV}), and all its groups $E^1_{p,q}$ consist of homology
groups of certain complex flag manifolds.

We consider a different spectral sequence, also converging to the same group
(\ref{alex}), also built of homology groups of flag manifolds, and based on a
{\em conical resolution} of $\Sigma.$ This sequence (and the related filtration
in the ring (\ref{split})) seems to be interesting, because
\begin{enumerate}
\item it degenerates at the {\em first} term
$E^1 \equiv E^\infty$ (at least in characteristic 0);
\item it is compatible very much with
the inclusions $\HH(n) \hookrightarrow \HH(n+1) \hookrightarrow \cdots,$ thus
stabilizing to a similar spectral sequence, converging to the cohomology ring
of the space of infinite Hermitian matrices without multiple eigenvalues, all
whose terms $E^{p,q}_r$ are finitely generated;
\item it is (conjecturally) multiplicative;
\item it is closely related with the shift operator of the
spectrum, and the corresponding filtration in the cohomology is invariant under
this operator;
\item it allows us to apply the (more or less easy or well--known)
facts on the topology of flag manifolds to the problems of geometrical
combinatorics, especially concerning topological partially ordered sets and
order complexes;
\item its algebraic presentation is very similar to
the theory of finite-order invariants of knots (the stratum of operators with
eigenvalues of multiplicities $a_1, \ldots, a_l$ corresponds to that of smooth
maps $S^1 \to \R^3$ with $l$ selfintersection points of the same
multiplicities), thus inventing to the ring (\ref{split}) useful structures and
notions of this theory. In particular, it allows us to define the orders of
cohomology classes, and, for any class of order $p$, its {\em symbol} (or {\em
generalized residue}) at the strata of $\Sigma(n)$ of complexity $p$ in the
same way as it was done in \cite{bjorn}, \cite{phasis} for knot invariants and
strata of singular knots.
\end{enumerate}

\section{Hermitian matrices with simple spectra}
\label{hermss}

If all eigenvalues $\lambda_i$ of a Hermitian operator $\C^n \to \C^n$ are
different, then they (and the corresponding one-dimensional complex
eigenspaces) can be ordered by $\lambda_1 < \cdots < \lambda_n.$ These
eigenspaces form $n$ line bundles over the space $\HH(n)\sm \Sigma(n)$ of all
such operators; let $c^1, \ldots, c^n$ be the first Chern classes of these
bundles.

\prop (see \cite{borel}). There is canonical ring isomorphism
\begin{equation}
\label{bott} H^*(\HH(n)\sm \Sigma(n)) \simeq \Z[c^1, \ldots, c^n]/\{Sym\},
\end{equation}
where $\{Sym\}$ is the ideal generated by all symmetric polynomials of positive
degrees. \quad $\Box$ \eprop

It is convenient to consider all such rings for all $n$ simultaneously, i.e. to
consider the ring
\begin{equation}
\label{infin} \Z[[a^0,a^1,a^{-1},a^2,a^{-2},\ldots]]/\{Sym\}
\end{equation}
of formal power series of infinitely many (in both directions) two-dimensional
variables $a^j$ factored through the ideal spanned by symmetric series of
positive degrees. This ring is called the cohomology ring of the space of {\em
infinite} Hermitian matrices with simple spectra.

The ring (\ref{bott}) can be identified in many ways with a quotient ring of
(\ref{infin}): we can choose any number $i=0,1, -1, \ldots$ and factor
(\ref{infin}) additionally through all elements $a^j$ with $j\le i$ and
$j>i+n$. Identifying then the variables $c^1, \ldots, c^n$ with $a^{i+1},
\ldots, a^{i+n}$ respectively, we get an isomorphism between this quotient ring
and (\ref{bott}).
\medskip

The obvious operator, mapping any variable $a^i$ to $a^{i+1}$, acts on the
algebra (\ref{infin}); it is called the {\em shift operator.}

\section{Topological posets and conical resolutions}

Suppose that we have a stratified variety and wish to calculate its homology
groups. There are two main approaches to this problem. The obvious method of
{\em open strata} is as follows: we filter the variety by unions $S_i$ of
smooth strata of "complexity $\ge i$", and consider the corresponding spectral
sequence (whose term $E^1_{p,q}$ is the group $H_{p+q}(S_i,S_{i+1})$,
Poincar\'e dual to a cohomology group of the smooth manifold $S_i \sm
S_{i+1}$).

A different approach, modelling the combinatorial formula of inclusions and
exclusions, is as follows: first we consider a singularity resolution
$\pi:\tilde V \to V$ of entire variety $V$ (thus changing it over the singular
set), and then improve it over the {\em closures} of strata of smaller and
smaller dimensions in such a way that at the last step we get a space $V'$ with
a {\em proper} projection onto $V$ and contractible preimages of all points.
Then $V'$ is homotopy equivalent to the original space; in all important cases
this space $V'$ has a very transparent topological structure, in particular a
natural filtration, whose spectral sequence degenerates very fast, see e.g.
\cite{v91}, \cite{phasis}, \cite{novikov}.
\medskip

{\sc Example.} The group $H^2(\HH(n)\sm \Sigma(n))$ is $(n-1)$-dimensional and
consists of all sequences $(\alpha_1, \ldots, \alpha_n)$ of integer numbers
(i.e., of corresponding sums $\sum \alpha_i c^i$) factored through the constant
sequences.

The spectral sequence from \cite{hall} has on the corresponding line
$\{p+q=n^2-3 \}$ unique nontrivial term $E^{n-1,n^2-n-2}_\infty$, isomorphic to
$\Z^{n-1}$ and canonically generated by linking numbers with all $n-1$ smooth
strata of maximal dimension in $\Sigma(n)$, i.e. by "$\delta'$-like" sequences
of the form $(0, \ldots, 0,1,-1, 0, \ldots, 0)$. Our spectral sequence defined
below has $n-1$ one-dimensional groups on the same line, and the $p$-th term of
the corresponding filtration consists of all integer-valued polynomial
sequences of degree $\le p$ (modulo the constants). The stable filtration in
the 2-dimensional component of (\ref{infin}) also is finitely generated: it is
defined by polynomials of any degrees, in particular is invariant under the
shift operators.
\medskip

The explicit realization of our method is based on the notion of the {\em
topological order complex} and on the techniques of {\em conical resolutions}.
It generalizes the method of {\em simplicial resolutions}, see e.g.
\cite{viniti}, \cite{phasis}, applicable in the case, when all essential
singularities of the variety $V$ are finite (self-)intersections only.

\medskip
Instead of formal definitions, we present here the following illustration,
important for our further calculations.

\subsection{Determinant variety and homology of the group $U(n)$.}

Consider the space $Mat(\C^n)$ of all complex linear operators $\C^n \to \C^n;$
the variety in question is its {\em determinant subvariety} $Det $ consisting
of all degenerate operators. Consider all possible kernels of degenerate
operators, i.e. all Grassmann manifolds $G_1(\C^n), $ $\ldots, G_{n-1}(\C^n),
G_n(\C^n).$ The incidence of subspaces, corresponding to their points, makes
the disjoint union of these Grassmannians a {\em poset} ($:=$partially ordered
set) with unique maximal element $\{\C^n\} \subset G_n(\C^n) $. Then we take
the join $G_1(\C^n) * \ldots * G_n(\C^n)$, i.e., roughly speaking, the union of
all simplices (of different dimensions), whose vertices belong to different
Grassmannians. Such a simplex is called {\em coherent}, if all subspaces in
$\C^n$ corresponding to its vertices form a flag. Finally, the {\em topological
order complex} $\Theta(n)$ is the subset of our join, defined as the union of
all coherent simplices. It is contractible: indeed, it is a cone with the
vertex $\{\C^n\}$. Its base $\partial \Theta(n)$ is defined in a similar way as
the union of all coherent simplices of the join $G_1(\C^n)* \ldots *
G_{n-1}(\C^n).$

\prop \label{theta} (see \cite{v91}, \cite{phasis}). The space $\partial
\Theta(n)$ is PL-homeomorphic to $S^{n^2-2}$. \label{km} \quad $\Box$ \eprop

{\sc Remark.} Probably this assertion (and its generalizations) is assumed in
Remark 1.4 of the work \cite{BS}. I thank M.~M.~Kapranov for this reference.
\medskip

For any subspace $L \subset \C^n$ we define the cone $\Theta(L)$ as the union
of coherent simplices subordinate to $L$, i.e. such that all subspaces
corresponding to their vertices belong to $L$. In particular, $\Theta(n) \equiv
\Theta(\{\C^n\}).$ The subspace $\chi(L) \subset Mat(\C^n)$ is defined as the
space of all operators, whose kernels contain $L$.

The desired {\em conical resolution} $\Delta_n$ of $Det$ is a subset of
$Mat(\C^n) \times \Theta(n)$. Namely, for any subspace $L \subset \C^n$ of any
positive dimension we consider the space $ \chi(L) \times \Theta(L) \subset
Mat(\C^n) \times \Theta(n) $ and define $\Delta_n$ as the union of such spaces
over all possible $L.$ This space is naturally filtered: the term $F_i$ of this
filtration is the union of sets $\chi(L) \times \Theta(L)$ over all $L$  of
dimensions $ \le i.$

\prop The obvious map $\Delta_n \to Det$ (defined by the projection $Mat(\C^n)
\times \Theta(n) \to Mat(\C^n)$) induces a homotopy equivalence of one-point
compactifications of these spaces, in particular an isomorphism of their
Borel--Moore homology groups. \eprop

Now consider the spectral sequence, calculating these groups and generated by
our filtration $\{F_i\}$ in $\Delta_n.$ By definition,
\begin{equation}
\label{sps} E^1_{p,q} \simeq \bar H_{p+q}(F_p \sm F_{p-1})
\end{equation}
$F_p \sm F_{p-1}$ is the space of a fiber bundle, whose base is $G_p(\C^n)$ and
the fiber over the point $\{L\}$ is diffeomorphic to the direct product
$\C^{n(n-p)} \times (\Theta(L) \sm \partial \Theta(L)).$ By Proposition
\ref{km} the last factor of this product is homeomorphic to the open disc of
dimension $p^2-1$, thus the group (\ref{sps}) is isomorphic to
$H_{t}(G_p(\C^n)),$ $t=p+q-(2n^2-2np+p^2-1).$

This spectral sequence degenerates at the first step: $E^1 \equiv E^\infty.$

The Alexander dual cohomological spectral sequence, defined by the identity
$E^{p,q}_r \equiv E_{-p,2n^2-q-1}^r,$ converges to the cohomology group of the
complementary space $GL(n,\C) \sim U(n).$ Its degeneration gives us immediately
the {\em Miller splitting formula}
\begin{equation}
\label{miller} H^m(U(n)) \cong \bigoplus_{p=0}^n H^{m-p^2}(G_p(\C^n)).
\end{equation}

Similar splittings hold also for other classical Lie groups $O(n)$ and $Sp(n)$,
cf. \cite{v91}, \cite{phasis}.
\medskip

This cohomological spectral sequence is multiplicative. Indeed, the ring
$H^*(U(n))$ is an exterior algebra with $n$ generators $\alpha_1, \alpha_3,
\ldots, \alpha_{2n-1}$ of corresponding dimensions. Our spectral sequence
induces a filtration in this ring, whose $i$-th term coincides with the
subalgebra generated by all products of $\le i$ elements $\alpha_j.$

\section{Conical resolution of the space of Hermitian matrices
with multiple spectra} \label{constr}

Consider again the discriminant variety $\Sigma(n) \subset \HH(n).$ The
corresponding conical resolution, order complexes etc. are constructed as
follows.

Let $A=(a_1 \ge a_2 \ge \cdots \ge a_l)$ be a sequence of naturals with $a_l
\ge 2$ and $\sum a_i \le n.$ The {\em complexity} of the index $A$ is the
number $\sum_{i=1}^l (a_i -1);$ its {\it length} $\#A$ is the number of its
elements $a_i$ (denoted above by $l$), and {\em liberty} $\delta(A)$ is the
number $n- \sum a_i$. Also set $|A|=\sum a_i.$

Denote by $\Gamma_A(n)$ the space of {\em unordered} collections of $\#A$
pairwise Hermitian--or\-tho\-go\-nal subspaces in $\C^n$ of complex dimensions
$a_1, \ldots, a_{\#A}.$ If all numbers $a_i$ are different, then it can be
identified with the flag manifold ${\mathcal F}_{a_1, a_1+a_2, \ldots,
a_1+\cdots +a_{\# A}},$ and if some of $a_i$ coincide then with the quotient
space of this flag manifold under the (free) action of corresponding
permutation group $S(A)$. In any case, it is a complex compact manifold.

The disjoint union of all spaces $\Gamma_A(n)$ is a partially ordered set: a
point $\gamma \in \Gamma_A(n)$ is subordinate to the point $\gamma' \in
\Gamma_{A'}(n)$ if any of $\#A$ subspaces forming the point $\gamma$ belongs to
some of $\#A'$ subspaces forming $\gamma'.$ In this case we say also that the
points $\gamma$ and $\gamma'$ are {\em incident} to one another. This poset has
unique maximal element: the point $\{\C^n\} \subset \Gamma_{(n)}(n)$.

The corresponding topological order complex $\Xi(n)$ is defined as in the
previous section: we consider the join of spaces $\Gamma_A(n)$ over all
possible multiindices $A$ with $|A| \le n$, and take the union of all coherent
simplices in it, i.e. of simplices, all whose vertices are incident to one
another, thus forming a monotone sequence. It is obviously a cone with the
vertex at the point $\{\C^n\} \subset \Gamma_{(n)}(n)$.

To any point $\gamma \in \Gamma_A(n)$ there corresponds the subcomplex
$\Xi(\gamma) \subset \Xi(n):$ the union of all coherent simplices, all whose
vertices are subordinate to $\gamma.$ In particular $\Xi(n) \equiv
\Xi(\{\C^n\}).$

Also define the {\em link} $\partial \Xi(\gamma)$ as the union of all coherent
simplices forming $\Xi(\gamma),$ which {\em do not} contain the maximal vertex
$\{\gamma\}.$ The {\em open cone} $\breve \Xi(\gamma)$ is the difference
$\Xi(\gamma) \sm \partial \Xi(\gamma).$ Their homology groups are related by
the boundary isomorphism $\partial: \bar H_*(\breve \Xi(\gamma))
\stackrel{\sim}{\to} \tilde H_{*-1}(\partial \Xi(\gamma)).$

\thm \label{link} For any $n$, the reduced homology group $\tilde H_i(\partial
\Xi(n), \C)$ is trivial in all even dimensions if $n$ is even and in all odd
dimensions if $n$ is odd. \etheorem

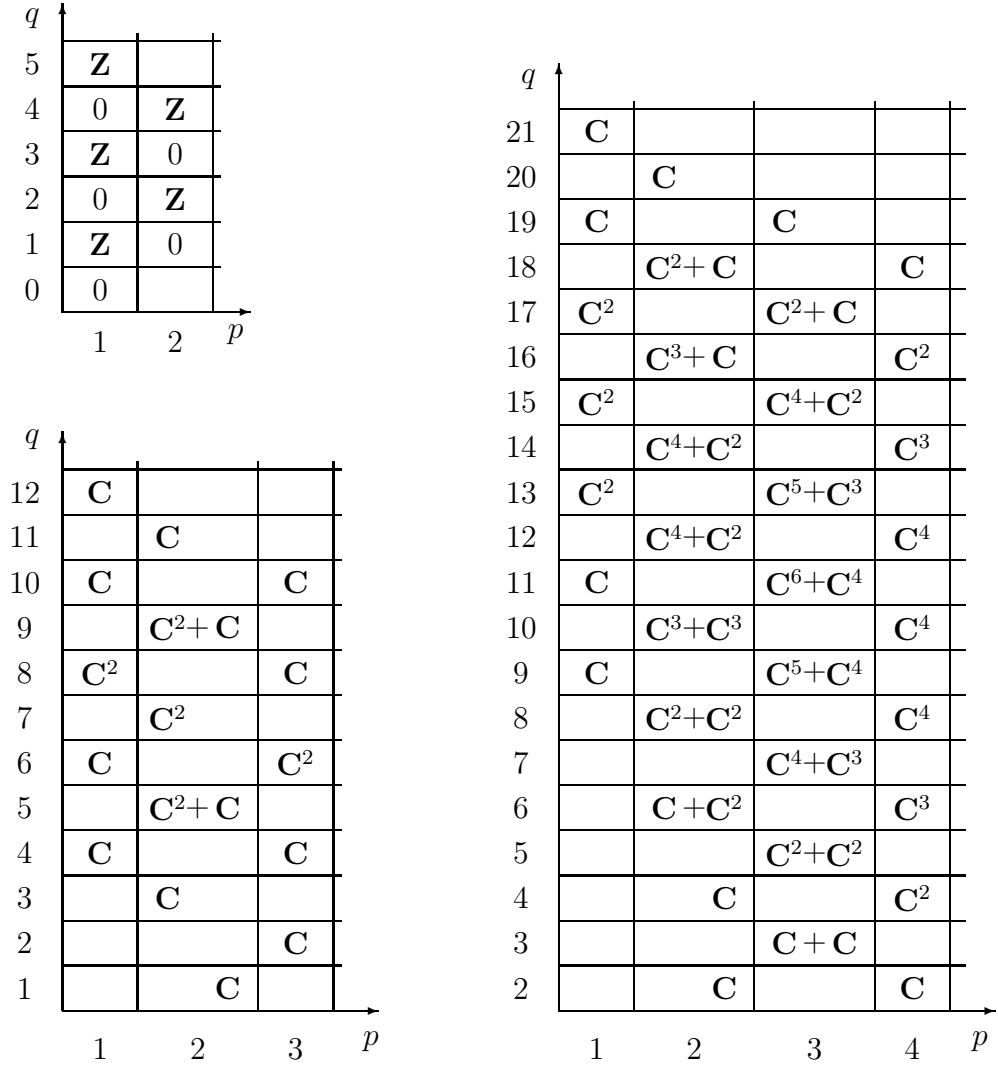
\begin{figure}
\unitlength=1.00mm \special{em:linewidth 0.4pt} \linethickness{0.4pt}
\begin{picture}(134.00,144.00)
\put(10.00,10.00){\vector(1,0){42.00}} \put(10.00,10.00){\vector(0,1){77.00}}
\put(5.00,13.00){\makebox(0,0)[cc]{1}} \put(5.00,19.00){\makebox(0,0)[cc]{2}}
\put(5.00,25.00){\makebox(0,0)[cc]{3}} \put(5.00,31.00){\makebox(0,0)[cc]{4}}
\put(5.00,37.00){\makebox(0,0)[cc]{5}} \put(5.00,43.00){\makebox(0,0)[cc]{6}}
\put(5.00,49.00){\makebox(0,0)[cc]{7}} \put(5.00,55.00){\makebox(0,0)[cc]{8}}
\put(5.00,61.00){\makebox(0,0)[cc]{9}} \put(5.00,67.00){\makebox(0,0)[cc]{10}}
\put(5.00,73.00){\makebox(0,0)[cc]{11}} \put(5.00,79.00){\makebox(0,0)[cc]{12}}
\put(20.00,10.00){\line(0,1){73.00}} \put(15.00,79.00){\makebox(0,0)[cc]{${\bf
C}$}} \put(15.00,67.00){\makebox(0,0)[cc]{${\bf C}$}}
\put(15.00,55.00){\makebox(0,0)[cc]{${\bf C}^2$}}
\put(15.00,43.00){\makebox(0,0)[cc]{${\bf C}$}}
\put(15.00,31.00){\makebox(0,0)[cc]{${\bf C}$}}
\put(24.00,25.00){\makebox(0,0)[cc]{${\bf C}$}}
\put(24.00,37.00){\makebox(0,0)[cc]{${\bf C}^2$}}
\put(24.00,49.00){\makebox(0,0)[cc]{${\bf C}^2$}}
\put(24.00,61.00){\makebox(0,0)[cc]{${\bf C}^2$}}
\put(24.00,73.00){\makebox(0,0)[cc]{${\bf C}$}}
\put(28.00,61.00){\makebox(0,0)[cc]{+}} \put(28.00,37.00){\makebox(0,0)[cc]{+}}
\put(32.00,13.00){\makebox(0,0)[cc]{${\bf C}$}}
\put(32.00,37.00){\makebox(0,0)[cc]{${\bf C}$}}
\put(32.00,61.00){\makebox(0,0)[cc]{${\bf C}$}}
\put(36.00,83.00){\line(0,-1){73.00}} \put(28.00,5.00){\makebox(0,0)[cc]{2}}
\put(15.00,5.00){\makebox(0,0)[cc]{1}}
\put(41.00,19.00){\makebox(0,0)[cc]{${\bf C}$}}
\put(41.00,31.00){\makebox(0,0)[cc]{${\bf C}$}}
\put(41.00,43.00){\makebox(0,0)[cc]{${\bf C}^2$}}
\put(41.00,55.00){\makebox(0,0)[cc]{${\bf C}$}}
\put(41.00,67.00){\makebox(0,0)[cc]{${\bf C}$}}
\put(46.00,83.00){\line(0,-1){73.00}} \put(41.00,5.00){\makebox(0,0)[cc]{3}}
\put(51.00,6.00){\makebox(0,0)[cc]{$p$}} \put(47.00,16.00){\line(-1,0){37.00}}
\put(10.00,22.00){\line(1,0){37.00}} \put(47.00,28.00){\line(-1,0){37.00}}
\put(10.00,34.00){\line(1,0){37.00}} \put(47.00,40.00){\line(-1,0){37.00}}
\put(10.00,46.00){\line(1,0){37.00}} \put(47.00,52.00){\line(-1,0){37.00}}
\put(10.00,58.00){\line(1,0){37.00}} \put(47.00,64.00){\line(-1,0){37.00}}
\put(10.00,70.00){\line(1,0){37.00}} \put(47.00,76.00){\line(-1,0){37.00}}
\put(10.00,82.00){\line(1,0){37.00}} \put(6.00,86.00){\makebox(0,0)[cc]{$q$}}
\put(10.00,103.00){\vector(1,0){25.00}} \put(10.00,103.00){\vector(0,1){41.00}}
\put(6.00,142.00){\makebox(0,0)[cc]{$q$}}
\put(25.00,130.00){\makebox(0,0)[cc]{${\bf Z}$}}
\put(25.00,124.00){\makebox(0,0)[cc]{0}}
\put(25.00,118.00){\makebox(0,0)[cc]{${\bf Z}$}}
\put(25.00,112.00){\makebox(0,0)[cc]{0}}
\put(15.00,106.00){\makebox(0,0)[cc]{0}}
\put(15.00,112.00){\makebox(0,0)[cc]{${\bf Z}$}}
\put(15.00,118.00){\makebox(0,0)[cc]{0}}
\put(15.00,124.00){\makebox(0,0)[cc]{${\bf Z}$}}
\put(15.00,130.00){\makebox(0,0)[cc]{0}}
\put(15.00,136.00){\makebox(0,0)[cc]{${\bf Z}$}}
\put(6.00,136.00){\makebox(0,0)[cc]{5}} \put(6.00,130.00){\makebox(0,0)[cc]{4}}
\put(6.00,124.00){\makebox(0,0)[cc]{3}} \put(6.00,118.00){\makebox(0,0)[cc]{2}}
\put(6.00,112.00){\makebox(0,0)[cc]{1}} \put(6.00,106.00){\makebox(0,0)[cc]{0}}
\put(15.00,99.00){\makebox(0,0)[cc]{1}} \put(25.00,99.00){\makebox(0,0)[cc]{2}}
\put(33.00,100.00){\makebox(0,0)[cc]{$p$}}
\put(10.00,109.00){\line(1,0){21.00}} \put(31.00,115.00){\line(-1,0){21.00}}
\put(10.00,121.00){\line(1,0){21.00}} \put(31.00,127.00){\line(-1,0){21.00}}
\put(10.00,133.00){\line(1,0){21.00}} \put(31.00,139.00){\line(-1,0){21.00}}
\put(20.00,140.00){\line(0,-1){37.00}} \put(30.00,103.00){\line(0,1){37.00}}
\put(76.00,10.00){\vector(1,0){58.00}} \put(76.00,10.00){\vector(0,1){126.00}}
\put(72.00,134.00){\makebox(0,0)[cc]{$q$}}
\put(71.00,13.00){\makebox(0,0)[cc]{2}} \put(71.00,19.00){\makebox(0,0)[cc]{3}}
\put(71.00,25.00){\makebox(0,0)[cc]{4}} \put(71.00,31.00){\makebox(0,0)[cc]{5}}
\put(71.00,37.00){\makebox(0,0)[cc]{6}} \put(71.00,43.00){\makebox(0,0)[cc]{7}}
\put(71.00,49.00){\makebox(0,0)[cc]{8}} \put(71.00,55.00){\makebox(0,0)[cc]{9}}
\put(71.00,61.00){\makebox(0,0)[cc]{10}}
\put(71.00,67.00){\makebox(0,0)[cc]{11}}
\put(71.00,73.00){\makebox(0,0)[cc]{12}}
\put(71.00,79.00){\makebox(0,0)[cc]{13}}
\put(71.00,85.00){\makebox(0,0)[cc]{14}}
\put(71.00,91.00){\makebox(0,0)[cc]{15}}
\put(71.00,97.00){\makebox(0,0)[cc]{16}}
\put(71.00,103.00){\makebox(0,0)[cc]{17}}
\put(71.00,109.00){\makebox(0,0)[cc]{18}}
\put(71.00,115.00){\makebox(0,0)[cc]{19}}
\put(71.00,121.00){\makebox(0,0)[cc]{20}} \put(76.00,16.00){\line(1,0){54.00}}
\put(130.00,22.00){\line(-1,0){54.00}} \put(76.00,28.00){\line(1,0){54.00}}
\put(130.00,34.00){\line(-1,0){54.00}} \put(76.00,40.00){\line(1,0){54.00}}
\put(130.00,46.00){\line(-1,0){54.00}} \put(76.00,52.00){\line(1,0){54.00}}
\put(130.00,58.00){\line(-1,0){54.00}} \put(76.00,64.00){\line(1,0){54.00}}
\put(130.00,70.00){\line(-1,0){54.00}} \put(76.00,76.00){\line(1,0){54.00}}
\put(130.00,82.00){\line(-1,0){54.00}} \put(76.00,88.00){\line(1,0){54.00}}
\put(130.00,94.00){\line(-1,0){54.00}} \put(76.00,100.00){\line(1,0){54.00}}
\put(130.00,106.00){\line(-1,0){54.00}} \put(130.00,112.00){\line(-1,0){54.00}}
\put(76.00,118.00){\line(1,0){54.00}} \put(130.00,124.00){\line(-1,0){54.00}}
\put(86.00,131.00){\line(0,-1){121.00}} \put(102.00,10.00){\line(0,1){121.00}}
\put(118.00,131.00){\line(0,-1){121.00}}
\put(132.00,6.00){\makebox(0,0)[cc]{$p$}}
\put(123.00,5.00){\makebox(0,0)[cc]{4}} \put(110.00,5.00){\makebox(0,0)[cc]{3}}
\put(94.00,5.00){\makebox(0,0)[cc]{2}} \put(81.00,5.00){\makebox(0,0)[cc]{1}}
\put(98.00,13.00){\makebox(0,0)[cc]{${\bf C}$}}
\put(98.00,25.00){\makebox(0,0)[cc]{${\bf C}$}}
\put(98.00,37.00){\makebox(0,0)[cc]{${\bf C}^2$}}
\put(98.00,49.00){\makebox(0,0)[cc]{${\bf C}^2$}}
\put(98.00,61.00){\makebox(0,0)[cc]{${\bf C}^3$}}
\put(98.00,73.00){\makebox(0,0)[cc]{${\bf C}^2$}}
\put(98.00,85.00){\makebox(0,0)[cc]{${\bf C}^2$}}
\put(98.00,97.00){\makebox(0,0)[cc]{${\bf C}$}}
\put(98.00,109.00){\makebox(0,0)[cc]{${\bf C}$}}
\put(71.00,127.00){\makebox(0,0)[cc]{21}}
\put(81.00,127.00){\makebox(0,0)[cc]{${\bf C}$}}
\put(76.00,130.00){\line(1,0){54.00}}
\put(90.00,121.00){\makebox(0,0)[cc]{${\bf C}$}}
\put(81.00,115.00){\makebox(0,0)[cc]{${\bf C}$}}
\put(90.00,109.00){\makebox(0,0)[cc]{${\bf C}^2$}}
\put(94.00,109.00){\makebox(0,0)[cc]{+}}
\put(106.00,115.00){\makebox(0,0)[cc]{${\bf C}$}}
\put(81.00,103.00){\makebox(0,0)[cc]{${\bf C}^2$}}
\put(81.00,91.00){\makebox(0,0)[cc]{${\bf C}^2$}}
\put(81.00,79.00){\makebox(0,0)[cc]{${\bf C}^2$}}
\put(81.00,67.00){\makebox(0,0)[cc]{${\bf C}$}}
\put(81.00,55.00){\makebox(0,0)[cc]{${\bf C}$}}
\put(90.00,37.00){\makebox(0,0)[cc]{${\bf C}$}}
\put(94.00,37.00){\makebox(0,0)[cc]{+}}
\put(90.00,49.00){\makebox(0,0)[cc]{${\bf C}^2$}}
\put(90.00,61.00){\makebox(0,0)[cc]{${\bf C}^3$}}
\put(90.00,73.00){\makebox(0,0)[cc]{${\bf C}^4$}}
\put(90.00,85.00){\makebox(0,0)[cc]{${\bf C}^4$}}
\put(90.00,97.00){\makebox(0,0)[cc]{${\bf C}^3$}}
\put(94.00,97.00){\makebox(0,0)[cc]{+}} \put(94.00,85.00){\makebox(0,0)[cc]{+}}
\put(94.00,73.00){\makebox(0,0)[cc]{+}} \put(94.00,61.00){\makebox(0,0)[cc]{+}}
\put(94.00,49.00){\makebox(0,0)[cc]{+}}
\put(106.00,19.00){\makebox(0,0)[cc]{${\bf C}$}}
\put(106.00,31.00){\makebox(0,0)[cc]{${\bf C}^2$}}
\put(106.00,43.00){\makebox(0,0)[cc]{${\bf C}^4$}}
\put(106.00,55.00){\makebox(0,0)[cc]{${\bf C}^5$}}
\put(106.00,67.00){\makebox(0,0)[cc]{${\bf C}^6$}}
\put(106.00,79.00){\makebox(0,0)[cc]{${\bf C}^5$}}
\put(106.00,91.00){\makebox(0,0)[cc]{${\bf C}^4$}}
\put(106.00,103.00){\makebox(0,0)[cc]{${\bf C}^2$}}
\put(114.00,103.00){\makebox(0,0)[cc]{${\bf C}$}}
\put(110.00,103.00){\makebox(0,0)[cc]{+}}
\put(110.00,91.00){\makebox(0,0)[cc]{+}}
\put(110.00,79.00){\makebox(0,0)[cc]{+}}
\put(110.00,67.00){\makebox(0,0)[cc]{+}}
\put(110.00,55.00){\makebox(0,0)[cc]{+}}
\put(110.00,43.00){\makebox(0,0)[cc]{+}}
\put(110.00,31.00){\makebox(0,0)[cc]{+}}
\put(110.00,19.00){\makebox(0,0)[cc]{+}}
\put(114.00,19.00){\makebox(0,0)[cc]{${\bf C}$}}
\put(114.00,31.00){\makebox(0,0)[cc]{${\bf C}^2$}}
\put(114.00,43.00){\makebox(0,0)[cc]{${\bf C}^3$}}
\put(114.00,55.00){\makebox(0,0)[cc]{${\bf C}^4$}}
\put(114.00,67.00){\makebox(0,0)[cc]{${\bf C}^4$}}
\put(114.00,79.00){\makebox(0,0)[cc]{${\bf C}^3$}}
\put(114.00,91.00){\makebox(0,0)[cc]{${\bf C}^2$}}
\put(128.00,10.00){\line(0,1){121.00}}
\put(123.00,13.00){\makebox(0,0)[cc]{${\bf C}$}}
\put(123.00,25.00){\makebox(0,0)[cc]{${\bf C}^2$}}
\put(123.00,37.00){\makebox(0,0)[cc]{${\bf C}^3$}}
\put(123.00,49.00){\makebox(0,0)[cc]{${\bf C}^4$}}
\put(123.00,61.00){\makebox(0,0)[cc]{${\bf C}^4$}}
\put(123.00,73.00){\makebox(0,0)[cc]{${\bf C}^4$}}
\put(123.00,85.00){\makebox(0,0)[cc]{${\bf C}^3$}}
\put(123.00,97.00){\makebox(0,0)[cc]{${\bf C}^2$}}
\put(123.00,109.00){\makebox(0,0)[cc]{${\bf C}$}}
\end{picture}
\caption{Main spectral sequences for $n=3$ (left top), $n=4$ (left bottom) and
$n=5$ (right)} \label{ss}
\end{figure}

{\sc Examples.} For $n=3, 4$ and $5,$ the Poincar\'e polynomials of groups
$\tilde H_*(\partial \Xi(n), \C)$ are equal to $t^2(1+t^2),$
$t^3(1+t^4)(1+t^2+t^4)$ and $t^4(1+t^2+t^4+t^6)(1+t^2+t^4+t^6+t^8+t^{10}),$
respectively. These groups are encoded also in the very right nontrivial
columns in three tables of Fig.~\ref{ss}.
\medskip

A recurrent method for calculating these homology groups will be described in
the end of \S~\ref{exfive}, however I do not know any compact expression for
them.

For any index $A$ and any $\gamma \in \Gamma_A(n)$ denote by $\chi(\gamma)$ the
subspace in $\HH(n)$ consisting of all Hermitian operators such that any of
spaces of dimension $a_i,$ forming the point $\gamma,$ belongs to an eigenspace
of this matrix. This is a real vector space of dimension $\#A+ (\delta(A))^2.$

The conical resolution $\sigma(n)$ of the discriminant variety $\Sigma(n)
\subset \HH(n)$ is a subset of the direct product $\HH(n) \times \Xi(n),$
namely, the union of products $\chi(\gamma) \times \Xi(\gamma)$ over all
possible indices $A$ and all points $\gamma \in \Gamma_A(n)$.

\prop The space $\sigma(n)$ is semialgebraic. The obvious map $\sigma(n) \to
\Sigma(n)$ (defined by the projection $\HH(n) \times \Xi(n) \to \HH(n)$) is
proper and induces a homotopy equivalence of one-point compactifications of
these spaces, in particular an isomorphism of their Borel--Moore homology
groups. \quad $\Box$ \eprop

There is {\em standard filtration} $\sigma_1(n) \subset \cdots \subset
\sigma_{n-1}(n) \equiv \sigma(n)$ in the space $\sigma(n)$: its term
$\sigma_i(n)$ is the union of products $\chi(\gamma) \times \Xi(\gamma)$ over
all $\gamma \in \Gamma_A(n),$ where $A$ is some multiindex of complexity $\le
i.$

The {\em main spectral sequence,} calculating the Borel--Moore homology group
of $\sigma(n)$, is that generated by this filtration; by definition its term
$E^1_{p,q}$ is isomorphic to $\bar H_{p+q}(\sigma_p(n) \sm \sigma_{p-1}(n)).$

\thm \label{degen} The main spectral sequence, calculating the group
(\ref{alex}) with complex coefficients, degenerates at the first term: $E^1
\equiv E^\infty.$ \etheorem

These spectral sequences for $n=3,4$ and $5$ are shown in Fig.~\ref{ss}.

Over the integers, the similar statement is not true: see the end of
\S~\ref{exfour}.

\medskip
Now, let us describe the terms $E^1_{p,q}$ of this spectral sequence. For any
point $\gamma = (\gamma_1, \ldots, \gamma_{\#A}) \in \Gamma_A(n)$ denote by
$\gamma^\perp$ the Hermitian-orthogonal complement of the linear span of all
subspaces $\gamma_i.$

\prop \label{bundl} Any space $\sigma_p(n) \sm \sigma_{p-1}(n)$ is a disjoint
union (over all indices $A$ of complexity exactly $p$) of spaces of fiber
bundles, whose bases are the spaces $\Gamma_A(n)$ and the fiber over the point
$\gamma \in \Gamma_A(n)$ splits into the direct sum of three spaces
\begin{equation}
\R^{\#A} \times \HH(\delta(A)) \times \breve \Xi(\gamma). \label{fiber}
\end{equation}
\eprop

Namely, the factor $\R^{\#A}$ is defined by eigenvalues of operators $\Lambda
\subset \chi(\gamma)$ on all subspaces $\gamma_i$ forming the point $\gamma;$
the factor $\HH(\delta(A))$ is defined by restrictions of such operators to the
orthogonal subspaces $\gamma^\perp$, and the sign $\, \breve{}$ over
$\Xi(\gamma)$ is due to the fact that the subspace $\chi(\gamma) \times
\partial \Xi(\gamma)$ lies in the smaller term $\sigma_{p-1}(n)$ of the
filtration.

Denote the spaces of these fiber bundles by $\beta_A(n)$. So, $\beta_A(n)$ is
the space of a fibered product of three bundles over $\Gamma_A(n)$, whose
fibers are three factors in (\ref{fiber}).

\prop The second bundle over $\Gamma_A(n)$, formed by spaces isomorphic to
$\HH(\delta(A))$, is orientable. \eprop

Indeed, it is induced from a similar bundle over the simply-connected
Grassmannian manifold $G_{|A|}(\C^n)$ by the obvious map (sending any
collection $\gamma$ of subspaces to their linear span).

On the other hand, the bundle of first factors in (\ref{fiber}) is
nonorientable if some of numbers $a_i$ coincide. Any such factor splits into
the sum of 1-dimensional subspaces associated with all $\#A$ subspaces forming
the basepoint $\gamma;$ a path in $\Gamma_A(n)$ providing an odd permutation of
these subspaces violates the orientation of this factor.

The following proposition reduces the homology of the space $\breve \Xi(A)$ to
these for indices $A'$ of length 1.

Recall the rigorous definition of the join $X * Y$ of two topological spaces
$X,Y$ as the quotient space of the product $X \times [-1,1] \times Y$ by the
equivalence relations $(x,-1,y) \sim (x,-1,y')$ and $(x,1,y) \sim (x',1,y)$ for
any $x,x' \in X$ and any $y,y' \in Y$.

\prop \label{union} Suppose that $A=A' \cup A'',$ $\gamma \in \Gamma_A(n),$ and
the collection of subspaces in $\C^n$ defining the point $\gamma$ splits into
the union of similar collections defining some points $\gamma' \in
\Gamma_{A'}(n) $ and $\gamma'' \in \Gamma_{A''}(n).$ Then there is a standard
homeomorphism of the complex $\Xi(A)$ to the join $ \Xi(\gamma') *
\Xi(\gamma')$, sending the vertex $\{\chi(\gamma)\}$ to the point
$\{\chi(\gamma')\} \times 0 \times \{\chi(\gamma'')\}$ and the link $\partial
\Xi(\gamma)$ to the image (under the above-described factorization map) of the
boundary $ (\Xi(\gamma') \times [-1,1] \times \partial \Xi(\gamma'')) \cup
(\partial \Xi(\gamma') \times [-1,1] \times \Xi(\gamma''))$ of the prism $
\Xi(\gamma') \times [-1,1] \times \Xi(\gamma'')$. \eprop

The proof is immediate. \quad $\Box$
\medskip

{\sc Corollary.} {\em Define the graded group $h_*(\gamma)$ by the equation
$h_i(\gamma) \equiv \bar H_{i-1}(\breve \Xi(\gamma), \C)$. Then for any
$A=(a_1, \ldots, a_{\# A})$ and $\gamma \in \Gamma_A(n)$ there is {\em almost}
canonical isomorphism
\begin{equation}
\label{tens} h_*(\gamma) \simeq \bigotimes_{i=1}^{\#A} h_*(\gamma_i),
\end{equation}
where all $\gamma_i$ are subspaces of dimensions $a_i$ forming the point
$\gamma$ and considered as the points of corresponding Grassmannian manifolds
$\Gamma_{(a_i)}(n).$}
\medskip

"Almost" here means the following. The isomorphism (\ref{tens}) depends on the
order of subspaces $\gamma_i$ forming the collection $\gamma.$ A reordering of
such subspaces (of equal dimensions $a_i$) takes this isomorphism to itself
multiplied by $(-1)^s,$ where $s$ is the parity of this permutation. In
particular, we cannot realize the isomorphism (\ref{tens}) by a natural
construction depending continuously on the point $\gamma \in \Gamma_A(n):$ a
path in $\Gamma_A(n)$ providing such an odd permutation will violate this
construction.
\medskip

Besides this bundle $\beta_A(n) \to \Gamma_A(n)$, it is convenient to consider
also its "universal covering"
\begin{equation}
\label{unicov} \beta!_A(n) \to \Gamma!_A(n).
\end{equation}

Here $\Gamma!_A(n)$ is the universal covering space of $\Gamma_A(n).$ This is a
principal covering, whose group $S(A)$ is the direct product of permutation
groups $S_\mu$, where $\mu$ are multiplicities, with which equal numbers $a_i$
meet in the sequence $A$. The space $\Gamma!_A(n)$ is diffeomorphic to the flag
manifold $\F_{a_1, a_1+a_2, \ldots, |A|}$. The bundle (\ref{unicov}) is again
the fibered product of three bundles with fibers indicated in (\ref{fiber}); in
this case the first bundle with fibers $\R^{\#A}$ is trivial, and the
decomposition (\ref{tens}) of the homology of the third bundle can be realized
so that it depends continuously on the basepoint $\gamma! \in \Gamma!_A(n).$
The group of the covering $\Gamma!_A(n) \to \Gamma_A(n)$ acts naturally on the
space $\beta!_A(n),$ permuting corresponding fibers, and $\beta_A(n)$ can be
considered as the quotient space under this action.

\section{Examples}

In this section we calculate our spectral sequences for $n \le 5$, and
determine our filtration in the groups $H^2 (\HH(n) \sm \Sigma(n))$.

\subsection{The marginal columns of the spectral sequence.}
\label{margin}

For any $n,$ the first term $\sigma_1(n) $ of $\sigma(n)$ is the {\em canonical
resolution} of $\Sigma (n)$, i.e. the space of pairs \{a 2-plane in $\C^n$; an
operator whose restriction on this plane is scalar\}.

Its Borel--Moore homology group coincides with the usual homology group of the
Grassmannian manifold $G_2(\C^n)$ with grading shifted by $(n-2)^2+1,$ see the
columns $\{p=1\}$ of all three tables in Fig.~\ref{ss}.

Now, let us consider the very right columns $\{p=n-1\}.$ The unique index $A$
of complexity $n-1$ is equal to $(n).$ The space $\Gamma_A(n)$ in this case
consists of unique point $\{\C^n\}$, thus by Proposition \ref{bundl} we have
the isomorphism of homology groups,
\begin{equation}
\label{rightcol} \bar H_i(\sigma_{n-1}(n) \sm \sigma_{n-2}(n)) \equiv \bar
H_{i}(\breve \Xi(n) \times \R^1) \simeq \tilde H_{i-2}(\partial \Xi(n)).
\end{equation}

\subsection{Cases $n=2$ and $3$.}

If $n=2,$ then $\Sigma(n)$ consists of scalar matrices ${{\lambda~0}\choose
{0~\lambda}}$, $\lambda \in \R,$ thus $\HH(2)\sm \Sigma(2) \sim S^2$, see
\cite{hall}. The ingredients of the construction of \S~\ref{constr} are as
follows. The unique admissible index $A$ is equal to $(2),$ the corresponding
space $\Gamma_A(2)$ is a point $O$, the space $\Xi(O) \equiv \breve \Xi(O)$
also is a point, and $\chi(O)$ is a real line. Thus the (homological) spectral
sequence consists of unique element $E_{1,0} \simeq \Z.$

Now let be $n=3.$ In this case both columns are described in \S~\ref{margin}.
Namely, the column $\{p=1\}$ again contains the homology groups of the
canonical resolution $\sigma_1(n) \to \Sigma(n)$. This resolution is a local
diffeomorphism over all points of $\Sigma(n)$ other than the scalar operators,
and such operators are "blown up" to spaces $G_2(\C^3)
 \equiv \partial \Xi(3).$ The term $\sigma_2(n)
\sm \sigma_1(n)$ consists of open cones $ \breve C(G_2(\C^3)) \sim \breve
\Xi(3),$ spanning all these fibers $\partial \Xi(3).$ Thus we get the table in
the left upper part of Fig.~\ref{ss}.

\subsection{The calculations for $n=4$.}
\label{exfour}

For $n=4$ there are exactly two indices $A$ of complexity $2$, namely $(3)$ and
$(2,2).$ The corresponding blocks $\beta_A(4)$ are as follows.

For $A=(3),$ the space $\Gamma_3(n)$ is the Grassmannian manifold $G_3(\C^4)
\cong \CP^3,$ and for any point $\gamma$ of this space, the space $\partial
\Xi(\gamma)$ is equal to $\CP^2.$ Thus the spectral sequence of the fiber
bundle $\beta_{(3)}(4) \to \Gamma_{(3)}(4)$ is as follows: its term
$\E^2_{a,b}$ is isomorphic to $H_a(\CP^3, \tilde H_{b-3}(\CP^2)),$ where the
number $3$ in the lower index $b-3$ is $\#A+ (\delta(A))^2+$ (the loss of
dimensions in the boundary homomorphism). This spectral sequence obviously
degenerates in this term $\E^2$ and gives us the direct summand of the column
$\{p=2\}$ written in its left part.

For $A=(2,2),$ the space $\Gamma_A(4)$ is the quotient space of the
Grassmannian manifold $G_2(\C^4)$ under the involution sending any 2-plane into
its Hermitian--orthogonal plane. It is easy to calculate that the complex
homology group of this manifold is isomorphic to $\C$ in dimensions $0, 4$ and
$8$ and is trivial in other dimensions. Three factors of the fiber
(\ref{fiber}) of the bundle $\beta_A(4) \to \Gamma_A(4)$ in this case are equal
respectively to $\R^2,$ one point, and the open interval. (For any $\gamma \in
\Gamma_{(2,2)}(4)$ this intervals consists of two segments joining the point
$\gamma$ with two subordinate points of $G_2(\C^4).$ Their endpoints lying in
$G_2(\C^4)$ are of smaller filtration and hence are removed.) The generator of
the group $\pi_1(\Gamma_A(4)) \sim \Z_2$ violates the orientations of both the
bundle of spaces $\R^2$ and that of open intervals, hence $\bar H_i(\beta_A(4))
\simeq H_{i-3}(\Gamma_A(4)),$ which gives us the right-hand part of column
$\{p=2\}.$

Finally, let us calculate the column $\{p=3\},$ i.e. the homology group of the
link $\partial \Xi(4).$

Consider the subset $\Omega$ of this link, swept out by all coherent segments
connecting points of manifolds $G_2(\C^4)$ and $G_3(\C^4)$.

\begin{lemma}
The group $\tilde H_a(\Omega)$ is isomorphic to $\Z$ in dimensions 7, 9, 11,
and is trivial in all other dimensions.
\end{lemma}

{\em Proof.} Consider the space $\Theta(4)$ studied in Proposition \ref{theta}
(and, accordingly to this proposition, homeomorphic to $S^{14}$). The space
$\Omega$ is canonically embedded in this sphere, and by the construction is
{\em Spanier--Whitehead dual} (see \cite{v91}) to the space $\CP^3$ (also
embedded in it). Hence their reduced homology groups are related by the
Alexander duality in $S^{14}$, and our lemma is proved. \quad $\Box$
\medskip

Further, the space $\partial \Xi(4) \sm \Omega$ is the space of a nonorientable
fiber bundle over the space $\Gamma_{(2,2)}(4)$, whose fibers are open
intervals.

The complex Borel--Moore homology group of this space $M$ participates in the
(splitting) Smith exact sequence of the double covering $G_2(\C^4) \to
\Gamma_{(2,2)}(4)$,
$$
\cdots \stackrel{0}{\to} \bar H_{i+1}(M,\C) \to H_i(G_2(\C^4),\C) \to
H_i(\Gamma_{(2,2)}(4),\C) \stackrel{0}{\to} \bar H_i(M,\C) \to \cdots
$$
and hence can easily be calculated: it is equal to $\C$ in dimensions 3, 5 and
7, and is trivial in all other dimensions.

Since $\bar H_*(M) \equiv \tilde H_*(\partial \Xi(4), \Omega),$ the homology
group of $\partial \Xi(4)$ follows from the exact sequence of this pair. The
result is encoded in the right-hand column of the lower left table in
Fig.~\ref{ss}.
\medskip

{\sc Remark.} The similar integer-valued spectral sequence does nod degenerate
at the term $E^1:$ the right-hand part of its column $\{p=2\}$  contains a
2-torsion, which disappears in the limit homology group.

\subsection{The case $n=5$.}
\label{exfive}

Here we prove that the term $E^1$ of the main spectral sequence calculating
$\bar H_*(\Sigma(5),\C)$ is as shown in the right-hand table of Fig.~\ref{ss}.
The column $\{p=1\}$ is already justified: it coincides with the homology group
of $G_2(\C^5)$ up to a shift of dimensions.

The left part of the column $\{p=2\}$ corresponds to the block
$\beta_{(3)}(5)$: its term on the level $q$ coincides with the
$(q-4)$-dimensional component of the graded group $H_*(G_3(\C^5)) \otimes
\tilde H_*(G_2(\C^3))$. (Here $4= \#A+(\delta(A))^2+1-p$.)

The right-hand part of the same column contains the homology of the space
$\beta_{(2,2)}(5)$. Its base space $\Gamma_{(2,2)}(5)$ is fibered over $\CP^4$
with fiber $\Gamma_{(2,2)}(4).$ The spectral sequence of this fibration
obviously degenerates at the second term, so that $H_*(\Gamma_{(2,2)}(5))
\simeq H_*(\CP^4) \otimes H_*(\Gamma_{(2,2)}(4)).$ The fiber (\ref{fiber}) of
the bundle $\beta_{(2,2)}(5) \to \Gamma_{2,2}(5)$ is the direct product of
$\R^3$ and an open interval, and this fibration is orientable. Therefore by the
Thom isomorphism we get the right-hand part of the column $\{p=2\}$.

The left part of the column $\{p=3\}$ contains the homology groups of the space
$\beta_{(4)}(5)$, which is fibered over $\Gamma_{(4)}(5) \equiv \CP^4$ with the
fiber $\R^2 \times \breve \Xi(4)$. The homology groups of the last space
$\breve \Xi(4)$ are calculated in the previous subsection (see the column
$\{p=3\}$ of the corresponding spectral sequence). It follows from this
calculation, that terms $\E^2_{a,b}$ of the spectral sequence of the fibration
$\beta_{(4)}(5) \to \CP^4$ can be nontrivial only if $a$ is even and $b$ is
odd. Therefore this sequence degenerates and gives us the left part of the
column $\{p=3\}$.

The right-hand part of the same column contains the homology of the space
$\beta_{(3,2)}(5).$ Its basespace $\Gamma_{(3,2)}(5)$ coincides with
$G_3(\C^5)$, and the three factors of the fiber (\ref{fiber}) over some its
point are respectively $\R^2,$ one point, and the open cone over the suspension
$\Sigma (\CP^2),$ cf. Proposition~\ref{union}. Therefore the column $\{p=3\}$
also is justified.

Finally, the shape of the last column $\{p=4\}$ follows from formula
(\ref{split}) (providing all groups $\bar H_i(\Sigma(n), \C)$),
Theorem~\ref{degen} (stating that any such group is the direct sum of groups
$E^1_{p,i-p}$), and the columns calculated previously (providing all other
elements of these sums).

A similar consideration allows us to calculate all groups $\partial \Xi(n)$ if
we know similar groups with smaller values of $n$, and also the homology groups
of all spaces $\beta_A(n)$ with $A$ of complexities $\le n-2$. Concerning the
calculation of these groups, see \S~\ref{proofs} below.

\subsection{2-dimensional cohomology classes
of $\HH(n) \sm \Sigma(n)$.} \label{twodim}

It is convenient to replace our homological spectral sequence $E^r_{p,q} \to
\bar H_{p+q}(\Sigma(n))$, defined in \S~\ref{constr}, by its "Alexander dual"
cohomological sequence, obtained from it by the formal inversion of indices,
\begin{equation}
\label{invert} E^{p,q}_r \equiv E^r_{-p,n^2-q-1}.
\end{equation}

This spectral sequence lies in the second quadrant in the wedge $\{p \le 0, p+q
\ge 0\}$ and converges exactly to the group $H^*(\HH(n)\sm \Sigma(n))$.

\prop For any $n$ the integer-valued cohomological spectral sequence has
exactly $n-1$ nontrivial terms $E_1^{p,q}$ with $p+q=2,$ namely $E^{-1,3},$
$E^{-2,4},$ $E^{-n+1,n+1},$ all of which are isomorphic to $\Z.$ The
corresponding filtration in the group $H^2(\HH(n) \sm \Sigma(n))$ is as
follows: its term $F_p$ consists of all sums $\sum_{i=1}^n \alpha_i c^i,$ where
the sequence $\{\alpha_i\}$ coincides with a polynomial of degree $\le p$
taking integer values at integer points. \eprop

So, the quotient group $E^{-p,p+2}$ is generated by the basic polynomial
sequence of degree $p,$ $\alpha_i = i(i-1)\cdots (i-p+1)/p!,$ or, if we are
interested only in $\C$-valued cohomology, then just by the monomial $i^p$.

\prop For any $n$ the group $E^{-1,5}$ or first-order elements of the group
$H^4(\HH(n) \sm \Sigma(n),\Z)$ is one-dimensional and is generated by the
series $\sum_{i=1}^n i (c^i)^2.$ \eprop

The proofs are immediate. \quad $\Box$

\section{Proof of Theorems \protect\ref{link},
\protect\ref{degen}} \label{proofs}

We shall prove these theorems by induction over $n$.

\begin{lemma}
\label{parity} Suppose that Theorem \ref{link} is true for all $n \le n_0,$ and
all elements $a_i$ of the multiindex $A =(a_1, \cdots, a_{\#A})$ do not exceed
$n_0;$ let $\beta_A(N)$ be the corresponding block of the resolution of the
discriminant space $\Sigma(N) \subset \HH(N),$ and $\beta!_A(N)$ its universal
covering space, see the end of \S~\ref{constr}. Then

a) the spectral sequence of the fiber bundle $\beta!_A(N) \to \Gamma!_A(N),$
converging to the group $\bar H_*(\beta!_A(N)),$ degenerates at the second
term;

b) the groups $\bar H_*(\beta!_A(N),\C)$ and $\bar H_*(\beta_A(N),\C)$ are
trivial in all even (respectively, odd) dimensions if $N$ is even
(respectively, odd).
\end{lemma}

{\em Proof.} $\beta!_A(N)$ is the space of a fiber bundle over the
simply-connected manifold $\Gamma!_A(N)$ (all whose odd-dimensional homology
groups are trivial), with the fiber (\ref{fiber}), all whose
$(N-even)$-dimensional Borel--Moore homology groups also are trivial by our
assumption. This implies statement a) of the Lemma, and statement b) follows
from it because the complex homology group of the base of a finite covering is
"not greater" than that of its space. \quad $\Box$

\prop \label{step} If Theorem \ref{link} is true for all $n < N,$ then Theorem
\ref{degen} is true for $n=N$, i.e. the spectral sequence calculating the group
$\bar H_*(\Sigma(N),\C)$ degenerates at the first term. \eprop

This proposition implies the assertion of Theorem \ref{link} for $n=N$ and thus
completes the step of induction. Indeed, by this proposition the group $\bar
H_i(\Sigma(N),\C)$ splits into the direct sum of groups $E^1_{p,i-p}$. By the
formulas (\ref{alex}), (\ref{split}) this group is trivial in all dimensions
comparable with $N$ mod 2. Hence also the group $E^1_{N-1,i-(N-1)} \equiv \bar
H_i(\beta_{(N)}(N),\C) \equiv \bar H_i(\R^1 \times \breve \Xi(N), \C) \cong
\tilde H_{i-2}(\partial \Xi(N),\C)$ is trivial for such $i$.
\medskip

In the proof of Proposition \ref{step} we shall use the following version of
the Poincar\'e duality in the singular semialgebraic sets like
$\overline{\beta_A(N)},$ $\overline{\beta!_A(N)},$ and $\overline{\sigma(N)}$,
cf. \cite{GM}.

Given a compact Whitney stratified semialgebraic variety $V$, we embed it into
some sphere $S^M$ as an absolute retract of some its neighborhood $U \subset
S^M$ such that its boundary $\partial U$ is a smooth manifold. Then a {\em
quasicycle} in $V$ is any relative cycle of the pair $(U, \partial U),$ generic
with respect to the stratification of $V$. Any complex-valued cohomology class
of the space $V$ can be realized as the intersection index with some such
quasicycle. The {\em support} of a quasicycle is its intersection set with $V$.
\medskip

{\em Proof of Proposition \ref{step}.} We shall prove the dual (and thus
equivalent) statement concerning the cohomological spectral sequence
calculating the Borel--Moore {\em cohomology} of $\Sigma(N)$ (i.e. the reduced
cohomology of its one-point compactification $\overline{\Sigma(N)}$). Its term
$e^{p,q}_1$ is equal to $\tilde H^{p+q} \left(\overline{(\sigma_p(N) \sm
\sigma_{p-1}(N))},\C \right) \equiv \oplus_{|A|-\#A=p}\tilde
H^{p+q}\left(\overline{\beta_A(N)},\C \right),$ summation over all indices $A$
of complexity $p.$ Thus we need to prove that for any element $\omega$ of such
a group $\tilde H^{p+q} \left(\overline{\beta_A(N)},\C \right),$ all its
differentials $d^i(\omega) \in e^{p+i,q-i+1}$ are trivial.

Let $\omega!$ be the lifting of $\omega$ to the cohomology of $\beta!_A(N),$
and $[\omega!]$ the compact quasicycle in $\beta!_A(N)$ Poincar\'e dual to this
class. Its projection $[\omega] \subset \beta_A(N)$ is Poincar\'e dual to some
nonzero integer multiple of the class $\omega.$

Consider the open subset $reg \beta!_A(N) \subset \beta!_A(N),$ consisting of
all pairs of the form \{a Hermitian operator $\Lambda$; a point $\xi \in \breve
\Xi(\gamma!)$\} (where $\gamma! \in \Gamma!_A(N)$) such that the restriction of
$\Lambda$ onto the orthogonal plane $(\gamma!)^\perp \subset \C^N$ is an
operator with simple (i.e. $\delta(A)$-element) spectrum.

\begin{lemma}
For any class $\omega! \in \bar H^*(\beta!_A(N),\C),$ the Poincar\'e dual
quasicycle $[\omega!]$ can be chosen so that it support lies in $reg
\beta!_A(N).$
\end{lemma}

{\em Proof.} This assertion is equivalent to the following one: the
homomorphism
\begin{equation}
\label{reg} \bar H^*(\beta!_A(N), \C) \to \bar H^*(reg \beta!_A(N), \C),
\end{equation}
induced by the identical embedding, is monomorphic. To prove it, denote by
$\Upsilon_A(N)$ the total space of the subbundle of the fiber bundle
$\beta!_A(N) \to \Gamma!_A(N)$, whose fibers are the products of only first and
third factors in (\ref{fiber}). Then we have the commutative diagram of fiber
bundles
\begin{equation}
\label{comdiag}
\begin{array}{ccc}
\beta!_A(N) & \gets & reg \, \beta!_A(N) \\
\downarrow & & \downarrow \\
\Upsilon_A(N) & \stackrel{Id}{\gets} & \Upsilon_A(N)
\end{array}
\end{equation}
with fibers $\HH((\gamma!)^\perp)$ and $\HH((\gamma!)^\perp) \sm
\Sigma((\gamma!)^\perp)$ respectively.

The left bundle is an orientable $(\delta(A))^2$-dimensional vector bundle,
hence the second term of the corresponding spectral sequence, calculating the
Borel--Moore cohomology of $\beta!_A(N)$, has unique nontrivial row
$q=(\delta(A))^2,$ which coincides with the graded group $\bar
H^*(\Upsilon_A(N), \C)$; in particular this sequence degenerates at this term.
By the Corollary of Proposition \ref{union} and Theorem 1 (which we assume to
be true for all values of $n$ equal to elements $a_i$ of the multiindex $A$)
this graded group trivial in dimensions comparable with $|A|$ mod 2.

The similar spectral sequence for the right-hand bundle in (\ref{comdiag}) also
degenerates at the second term: this again follows from dimensional reasons,
because the group $\bar H^*(\HH(\gamma!^\perp) \sm \Sigma(\gamma!^\perp),\C)$
is nontrivial only in dimensions, comparable with $\delta(A)$ mod 2.

The homomorphism of these second terms of spectral sequences is an isomorphism
of rows $\{q=(\delta(A))^2\}$, and is zero operator for all other $q$. This
implies our lemma.
\medskip

Further, consider the subset $perf \beta!_A(N) \subset reg \beta!_A(N),$
consisting of pairs of the form \{an operator $\Lambda$; a point $\xi \in
\breve \Xi(\gamma!)$\}, where $\gamma!=(\gamma_1, \ldots, \gamma_{\#A})$, such
that additionally the eigenvalues of $\Lambda$ on the space $\gamma_1$
(respectively, $\gamma_2,$ \ldots, $\gamma_{\#A}$, $\gamma!^\perp$) lie in the
interval $(0,1)$ (respectively, $(2,3), \ldots, (2\#A-2, 2\#A-1), (2\#A,
+\infty)$).

\begin{lemma}
For any class $\omega! \in \bar H^*(\beta!_A(N),\C),$ the Poincar\'e dual
quasicycle $[\omega!]$ can be chosen so that it support lies in perf
$\beta!_A(N).$
\end{lemma}

{\em Proof.} On the space $\beta!_A(N)$ there acts the group $\R^{\#A} \oplus
\R,$ where $\R^{\#A}$ consists of the shifts along the first (trivial) factor
in (\ref{fiber}), and $\R$ of addings the operators which are scalar on the
spaces $\gamma!^\perp$ and trivial on all $\gamma_i.$ (If $\delta(A)=0$, then
the last action of $\R$ is trivial.) This action preserves all fibers of the
bundle $\beta!_A(N) \to \Gamma!_A(N)$ and the subspace $reg \beta!_A(N)$, and
takes quasicycles to quasicycles.

Any compact subset in $\beta!_A(N)$ obviously can be moved to the domain $perf
\beta!_A(N)$ by the action of this group and the additional action of the group
$\R_+$ of dilations in $\C^n.$ \quad $\Box$.
\medskip

So, our class $\omega \in H^*(\beta_A(N),\C)$ can be realized by the projection
$[\omega]$ of some compact quasicycle $[\omega!] \subset perf \beta!_A(N).$ But
the projection of the set perf $\beta!_A(N)$ to $\beta_A(N)$ does not meet the
closures of any other blocks $\beta_{A'}(N)$ of greater filtration (i.e., with
$|A'|-\#A'> |A|-\#A$). Thus our quasicycle $[\omega]$ is a quasicycle in entire
space $\overline{\sigma(n)}$, the dual Borel--Moore cohomology class in
$\sigma(n)$ is well defined, and Proposition \ref{step} is proved. \quad $\Box$

\section{Stabilization}

{\sc Definition.} The {\em order} of an element of the cohomology group
(\ref{alex}) is the filtration of the Alexander dual homology class in
$\Sigma(n)$, i.e. the smallest number $i$ such that it can be realized by a
locally finite chain inside $\sigma_i$.
\medskip

Consider again the cohomological spectral sequences $E_r^{p,q}(n) \to
H^{p+q}(\HH(n) \sm \Sigma(n)),$ see (\ref{invert}).

Let $n < N$ be any two naturals, $s$ one of numbers $0, \ldots, N-n,$ and $i_s:
\HH(n)\to \HH(N)$ an embedding, sending an operator $\Lambda : \C^n \to \C^n$
to the operator $\Lambda+\Lambda',$ where $\Lambda'$ acts on the orthogonal
complement of $\C^n$ in $\C^N,$ depends smoothly on $\Lambda$ and has a simple
spectrum, exactly $s$ (respectively, $N-n-s$) elements of which are below
(respectively, above) the spectrum of $\Lambda.$ In particular $\Sigma(n) =
i_s^{-1}(\Sigma(N))$.

\prop \label{stab} Any such inclusion induces a homomorphism of cohomological
spectral sequences,
\begin{equation}
\label{epi} E_r^{p,q}(N) \to E_r^{p,q}(n),
\end{equation}
$r \ge 1$ which {\em does not depend} on $s$ and on the inclusion $i_s$, is
epimorphic for any $p$ and $q$, and provides an isomorphism of terms
$E_r^{p,q}$ if $n$ is sufficiently large with respect to $|p|+|q|.$ \eprop

Thus we can define the {\em stable cohomological spectral sequence}, whose term
$E_r^{p,q}$ is the inductive limit of $E_r^{p,q}(n)$ over all possible such
homomorphisms. By the last statement of the proposition any such stable term is
finitely generated. This sequence converges to entire ring (\ref{infin}), thus
defining the notion of the order also for elements of this ring. In easy terms,
an element of the ring (\ref{infin}) is {\em of order} $i$ if for any
sufficiently large $n$ and any projection of (\ref{infin}) onto $H^*(\HH(n)\sm
\Sigma(n),\C)$, described in \S~\ref{hermss}, it becomes an element of order
$i.$ Certainly, not all elements of the ring (\ref{infin}) have some finite
order, however the spaces of finite-order elements {\em weakly converge} to
this ring in the following precise sense: for any $n$ any cohomology class on
the space $\HH(n)\sm \Sigma(n)$ coincides with the restriction of some stable
class of finite order. For the similar theory of knot invariants the similar
statement (the completeness of finite-order invariants) is not proved.
\medskip

Here is an estimate for the convergence of groups $E_r^{p,q}$. For any
multiindex $A$, the obvious inclusions of flag manifolds, $\Gamma!_A(n)
\hookrightarrow \Gamma!_A(n+1) \hookrightarrow \cdots$, induce the maps of
their cohomology groups. For any $i$ denote by $stab(A,i)$ the smallest number
$n$ at which all these cohomology groups of dimensions $\le i$ stabilize, i.e.
all further inclusions induce isomorphisms in these dimensions.

\prop \label{stabound} For any biindex $(p,q)$ with $p\le 0, p+q \ge 0$, the
groups $E_1^{p,q}(n)$ stabilize no later than at the instant $n =
\max_{|A|-\#A=-p} stab(A, p+q-2\#A).$ \eprop

The precise construction of the homomorphism mentioned in Proposition
\ref{stab} is as follows (cf. \cite{book}, \cite{phasis}). The embedding $i_s:
(\HH(n), \Sigma(n)) \to (\HH(N), \Sigma(N))$ can be tautologically lifted to
the filtration-preserving map of resolutions, $i_s: \sigma(n) \to \sigma(N).$
Its image admits in $\sigma(N)$ an open neighborhood $U$ homeomorphic to
$\sigma(n) \times \R^{N^2-n^2},$ where the image of any space $* \times
\R^{N^2-n^2}$ lies completely in one and the same block $\beta_A(N).$ Then we
get the map $\bar H_*(\sigma(N)) \to \bar H_{*-(N^2-n^2)}(\sigma(n))$ as the
composition of the restriction homomorphism $\bar H_*(\sigma(N)) \to \bar
H_{*}(U)$ and the K\"unneth isomorphism $\bar H_*(\sigma(n) \times
\R^{N^2-n^2}) \to \bar H_{*-(N^2-n^2)}(\sigma(n))$. By the construction, it is
compatible with the map of Alexander dual groups, $i_s^*: H^*(\HH(N)\sm
\Sigma(N)) \to H^*(\HH(n) \sm \Sigma(n))$ induced by the embedding $i_s.$ The
homomorphism of terms $E_1$ of our spectral sequences is defined by the
restrictions of this construction on any block $\beta_A(N)$.

For any $A$, set $U_A \equiv U \cap \beta_A(N)$. The homeomorphism $\beta_A(n)
\times \R^{N^2-n^2} \to U_A$ is compatible with the structure of the fiber
bundle in spaces $\beta_A(\cdot)$, mentioned in Proposition \ref{bundl}.
Namely, the embedding $\C^n \hookrightarrow \C^N$ defines the embeddings
$\Gamma_A(n) \hookrightarrow \Gamma_A(N)$ and $\Gamma!_A(n) \hookrightarrow
\Gamma!_A(N)$; there is a tubular neighborhood $W$ of $\Gamma_A(n) $ in $
\Gamma_A(N)$, whose bundle is orientable and fibers are homeomorphic to
$\R^{2|A|(N-n)}.$ We get the diagrams of bundles
\begin{equation}
\label{diags}
\begin{array}{ccccc}
\beta!_A(n) & \hookrightarrow & U! & \subset & \beta!_A(N)\\
\downarrow & & \downarrow & & \downarrow \\
\Gamma!_A(n) & \hookrightarrow & W! & \subset & \Gamma!_A(N),
\end{array}
\hspace{2cm}
\begin{array}{ccccc}
\beta_A(n) & \hookrightarrow & U & \subset & \beta_A(N)\\
\downarrow & & \downarrow & & \downarrow \\
\Gamma_A(n) & \hookrightarrow & W & \subset & \Gamma_A(N).
\end{array}
\end{equation}
For any point $\gamma \in \Gamma_A(n)$ or $\gamma! \in \Gamma!_A(n)$, the
fibers $\R^{\#A}$ and $\breve \Xi(\gamma)$ of the first and third factors of
the fibered products $\beta_A(n) \to \Gamma_A(n)$, $\beta!_A(n) \to
\Gamma!_A(n)$ over this point (see (\ref{fiber})) go identically to the
corresponding fibers of similar bundles over $\Gamma_A(N)$ and $\Gamma!_A(N)$,
while the fibers $\HH(n-|A|)$ of the second (vector) bundle become embedded
into similar fibers $\HH(N-|A|)$ of the second bundle over $\Gamma!_A(N)$ as
subspaces, whose quotient spaces form an orientable
$(N^2-n^2-2|A|(N-n))$-dimensional vector bundle over the image of $\Gamma_A(n)$
or $\Gamma!_A(n)$. The permutation group $S(A)$ acts on the left-hand diagram
(\ref{diags}), and the right-hand diagram is formed by the spaces of its
orbits; thus all complex homology homomorphisms induced by the arrows in the
latter diagram are determined by similar homomorphisms for the former one.

But these homomorphisms in the complex Borel--Moore homology groups are very
easy: the bottom homomorphism
\begin{equation}
\label{poin} \bar H_{*}(\Gamma!_A(n), \C) \gets \bar
H_{*+2|A|(N-n)}(\Gamma!_A(N), \C)
\end{equation}
is nothing but the standard cohomology map, induced by the obvious embedding
and rewritten in the terms of Poincar\'e dual groups. The Borel--Moore homology
groups of spaces of both bundles of this diagram are just the tensor products
of corresponding groups (\ref{poin}) of their bases and the homology groups of
fibers, which coincide up to a shift of dimensions by $N^2-n^2-2|A|(N-n).$ The
corresponding map $\bar H_{*}(\beta!_A(n), \C) \gets \bar
H_{*+(N^2-n^2)}(\beta!_A(N), \C)$ of these tensor products decomposes into
these actions on their factors.

The map (\ref{poin}) is always epimorphic, thus also the map (\ref{epi}) is.
\medskip

{\em Proof of Proposition \ref{stabound}.} Let us fix some multiindex $A$ and
dimension $i$ and estimate the greatest dimension $J=J(A,i)$ such that the
group $\bar H_{J}(\beta!_A(n))$ depends on the group $H^i(\Gamma!_A(n)).$ By
Proposition \ref{bundl} and formula (\ref{tens}) this dimension $J$ does not
exceed the sum of 4 numbers: $\dim \Gamma!_A(n)-i$, $\#A$, $(n-|A|)^2$ and
$\#A-1+\sum_{i=1}^{\#A} \dim_h \breve \Xi(a_i),$ where $\dim_h(\cdot)$ is the
maximal dimension $i$ such that $\bar H_i(\cdot)$ is nontrivial.

Define $s_2(A)$ as the second symmetric polynomial of numbers $a_1, \ldots,
a_{\#A},$ $s_2(A)= a_1a_2 + a_1a_3 + \cdots + a_{\#A-1}a_{\#A},$ then $\dim
\Gamma!_A(n) = 2(s_2(A)+ |A|(n-|A|))$. By \S~\ref{margin} the group $\bar
H_i(\breve \Xi(a),\C)$ is a direct summand of the group $\bar
H_{i+1}(\Sigma(a)).$ Since $\dim \Sigma(a) = a^2-3,$ we have $\dim_h \breve
\Xi(a) \le a^2-4.$ Thus the sum of four above numbers does not exceed
$2(s_2(A)+|A|(n-|A|))-i + \#A + (n-|A|)^2+\#A-1 + \sum a_i^2 - 4\#A \equiv
n^2-1-2\#A-i.$ Therefore by (\ref{invert}) any cell $E^{p,q}_1$ depends only on
such groups $H^i(\Gamma!_A(n),\C)$ that $|A|-\#A=-p$ and $p+q \ge i+ 2|A|.$
\quad $\Box$

\subsection*{Residues of cohomology classes at strata of
$\Sigma(n)$.}

{\sc Definition.} Given a class $\omega \in H^*(\HH(n)\sm \Sigma(n))$ of order
$p$ and a multiindex $A$ of complexity $p,$ the corresponding {\em symbol}
$s(A,\omega) \in \bar H_*(\beta_A(n))$ is the restriction of the corresponding
homology class in $\bar H_*(\Sigma)$ to the block $\beta_A(n)$. (In a more
detailed way, we realize the latter class by a cycle lying in $F_p(\sigma(n))$
and then reduce modulo the union of $\sigma_{p-1}$ and all other blocks
$\beta_{A'}(n)$ of complexity $p$, $A' \ne A$.) Similarly, the {\em !-symbol}
$s!(A,\omega) \in \bar H_*(\beta!_A(n))$ is the homology class of the complete
preimage of this class in the covering $\beta!_A(n)$.
\medskip

By Lemma \ref{parity}, any such !-symbol is the element of the cohomology group
of the flag manifold $\Gamma!_A(n)$ with coefficients in the homology group of
the topological order complex $\partial \Xi(A).$ In a similar way, the usual
symbol $s(A,\omega)$ is the element of the cohomology group of the manifold
$\Gamma_A(n)$ with coefficients in a local system of groups with the same
fibers.

\section{Problems}

{\bf Problem 1.} Our cohomological spectral sequence defines filtrations
("orders") in the rings (\ref{bott}), (\ref{infin}). The problem is to find the
explicit expression of elements of this filtration. I hope very much that it
coincides with or is closely related to some classical structure in these
rings.

In some sense, the term $F_{-p}$ of this filtration should consist of power
series, whose coefficients are "constructive functions of complexity $\le p$"
of corresponding exponents, cf. \S~\ref{twodim}. What is the precise sense of
this "complexity"?

Does there exist any natural (and physically motivated) presentation of
elements of this filtration in terms of differential forms in $\HH(n)\sm
\Sigma(n)$, cf. \cite{hall}?
\medskip

Here are two related subproblems, which can also be investigated independently.
\medskip

{\bf Problem 2.} The shift operator ${\bf S}$ (see the end of \S~\ref{hermss})
preserves our filtration: ${\bf S}(F_p) = F_p.$ Is it true, that the
corresponding derivative, sending the cohomology class $\alpha$ to ${\bf
S}(\alpha)-\alpha,$ maps any $F_p$ to $F_{p-1}$?
\medskip

{\bf Problem 3:} Is there true the following {\sc Multiplicativity conjecture:}
{\em Our cohomological spectral sequence is multiplicative, i.e.} $F_{p'}
\smile F_{p''} \subset F_{p'+p''}.$
\medskip

If it is true, the next problem is to define the corresponding maps
$E_1^{p',q'} \otimes E_1^{p'',q''} \to E_1^{p'+p'', q'+q''}$ explicitly.
\medskip

There is some obvious multiplication of this sort (in some sense modelling the
Kontsevich's formula for multiplication of knot invariants in the terms of
chord diagrams, see e.g. \cite{phasis}.) Namely, consider any multiindices
$A'=(a'_1, \ldots, a'_{\#A'})$ and $A''=(a''_1, \ldots, a''_{\#A''})$ and set
$A = A' \cup A''$. Below we define the natural map
\begin{equation}
\label{multipl} \bar H_{n^2-i-1} (\beta_{A'}(n),\C) \otimes \bar H_{n^2-j-1}
(\beta_{A''}(n),\C) \to \bar H_{n^2-i-j+1} (\beta_A(n),\C);
\end{equation}
the desired multiplication will follow from these by the Alexander duality
(\ref{invert}) and the decomposition $\sigma_p(n) \sm \sigma_{p-1}(n) =
\bigsqcup_{|A|-\#A=p} \beta_A(n).$

Indeed, define the group $\bar H_*(\breve \Xi(A))$ as $\bar H_*(\breve
\Xi(\gamma))$ for any $\gamma \in \Gamma_A(n)$. Then the map (\ref{multipl}) is
the composition of following maps of complex homology groups: \noindent
$$
\begin{array}{ccccccccc}
\bar H_*(\beta_{A'}(n)) & \stackrel{\triangle}{\to} & \bar H_*(\beta!_{A'}(n))
& \stackrel{\sim}{\to} & H_*(\Gamma_{A'}(n)) & &
\bar H_*(\R^{2\#A'-1+(n-|A'|)^2}) & & \bar H_*(\breve \Xi(A')) \\
\otimes & & \otimes & & \otimes & \otimes & \otimes & \otimes & \otimes\\
\bar H_*(\beta_{A''}(n)) & \stackrel{\triangle}{\to} & \bar
H_*(\beta!_{A''}(n)) & \stackrel{\sim}{\to} & H_*(\Gamma_{A''}(n)) & &
\bar H_*(\R^{2\#A''-1+(n-|A''|)^2}) & & \bar H_*(\breve \Xi(A'')) \\
& & & & \downarrow & & \downarrow & & \downarrow \\
\bar H_*(\beta_A(n)) & \stackrel{\nabla}{\gets} & \bar H_*(\beta!_A(n)) &
\stackrel{\sim}{\gets} & H_*(\Gamma_{A}(n)) & \otimes & \bar
H_*(\R^{2\#A-1+(n-|A|)^2}) & \otimes & \bar H_*(\breve \Xi(A)).
\end{array}
$$
Here the maps $\stackrel{\triangle}{\to}$ are the liftings, representing cycles
in $\beta_{A}(n)$ as projections of some cycles from $\beta!_{A}(n);$ all
identities $\stackrel{\sim}{\to},$ $\stackrel{\sim}{\gets}$ follow from
Proposition \ref{bundl} and statement a) of Lemma \ref{parity}. The vertical
map $ H_*(\Gamma!_{A'}(n)) \otimes H_*(\Gamma!_{A''}(n)) \to
H_*(\Gamma!_{A}(n))$ (attention! it does not preserves degrees!) is the
composition of a) Poincar\'e isomorphisms in $\Gamma!_{A'}(n)$,
$\Gamma!_{A''}(n)$, b) the K\"unneth isomorphism $H^*(\Gamma!_{A'}(n)) \otimes
H^*(\Gamma!_{A''}(n)) \to H^*(\Gamma!_{A'}(n) \times \Gamma!_{A''}(n))$, c) the
map in cohomology induced by the tautological embedding of $\Gamma!_{A}(n) $
into this product, and finally d) the Poincar\'e isomorphism in
$\Gamma_{A}(n)$. The vertical map $\bar H_*(\R^{2\#A'-1+(n-|A'|)^2}) \otimes
\bar H_*(\R^{2\#A''-1+(n-|A''|)^2}) \to \bar H_*(\R^{2\#A-1+(n-|A|)^2})$ (also
not preserving the grading) maps the product of canonical generators of two
first groups into the canonical generator of the last one. The vertical maps $
\bar H_*(\breve \Xi(A')) \otimes \bar H_*(\breve \Xi(A'')) \to \bar H_*(\breve
\Xi(A)) $ are described in the Proposition \ref{union} (they increase the
grading by 1), finally the map $\stackrel{\nabla}{\gets}$ is the obvious
projection. It is easy to check that this composition maps a
$(n^2-i+1)$-dimensional cycle in $\beta_{A'}(n)$ and a $(n^2-j+1)$-dimensional
cycle in $\beta_{A''}(n)$ into a $(n^2-i-j-1)$-dimensional homology class in
$\beta_A(n)$ and does not depend on the choice of liftings
$\stackrel{\triangle}{\to}$ of these cycles to corresponding spaces
$\beta!_{.}(n)$.
\medskip

{\sc Problem.} If the multiplication conjecture is true, how is the actual
multiplication in the spectral sequence related with the one just described?
What are its multiplicative generators?
\medskip

Besides our increasing filtration (by orders) the ring (\ref{infin}) admits a
decreasing filtration: its element is of {\em confinement} $d$ if for any $n<d$
it lies in the kernels of all surjections onto the ring (\ref{bott}), described
in \S~\ref{hermss}.

Is not it possible to characterize the primitive (i.e., indecomposable)
elements of our spectral sequence in the terms of these two filtrations?
\medskip

{\bf Problem 4.} To find the compact general formula for the homology of the
space $\partial \Xi(n)$ for any $n.$

This space is not a topological manifold (note that the right-hand columns of
tables in Fig.~\ref{ss} present the {\em reduced} homology of these spaces).

Therefore also their intersection homology groups are very interesting.

\subsection*{The last remark}

All results of this work can be immediately carried over to the theory of {\em
hyperhermitian} forms (see \cite{arrel}), i.e. of quadratic forms on the
realification $\R^{4n}$ of the quaternionic space ${\mathcal H}^n$ invariant
under the left multiplications by quaternions of length 1.

In particular, the main spectral sequence calculating the rational Borel--Moore
homology groups of the space of such forms with $<n$ eigenvalues degenerates at
the term $E^1$, and the reduced homology groups of the quaternionic analog of
the order complex $\partial \Xi(n)$ are trivial in dimensions not comparable
with $2n^2-n+1$ modulo 4.

I thank very much V.~I.~Arnold for conversations and remarks.

\end{document}